\documentclass [a4paper,english]{article}
\usepackage{amsthm, amsmath, amssymb, amsbsy, amsfonts, latexsym, euscript}
\usepackage[cp1251]{inputenc}
\usepackage[T2A]{fontenc}
\usepackage{indentfirst,epigraph}
\usepackage[english,russian]{babel}
\usepackage{graphicx,color}

\title{GINI MEANS AND ITS APPLICATIONS IN POLYMER CHEMISTRY}

\author{A.B. PEVNYI ,  S.M. SITNIK}

\date{}

\begin{document}
\maketitle

\selectlanguage{english}
\begin{abstract}

We consider Gini means with short biographical information and propose a new proof of the main inequality for these means. Also some applications of Gini and other means are considered to polymer chemistry.

\end{abstract}

\newpage

\selectlanguage{russian}

\begin{center}
{\Large \bf
СРЕДНИЕ ДЖИНИ И ИХ ПРИМЕНЕНИЯ В ХИМИИ ПОЛИМЕРОВ}
\end{center}

\begin{abstract}
Исследуются неравенства для средних, введённых итальянским математиком и статистиком Коррадо Джини. Приводится новое более простое доказательство основной теоремы для этих средних. Рассмотрены приложения средних Джини и некоторых других средних в теории полимеризации в химии. 
\end{abstract}

1. Широко известны три вида средних значений для пары положительных чисел $a,b$. Это \textsl{средние арифметическое} $A(a,b)$, \textsl{геометрическое} $G(a,b)$ и \textsl{гармоническое} $H(a,b)$, которые определяются по формулам
\begin{equation*}
A(a,b)=\frac{a+b}{2},\ \ \   G(a,b)=\sqrt{ab},\ \ \  H(a,b)=\frac{2ab}{a+b}.
\end{equation*}
Эти средние были введены ещё в Древней Греции в работах Пифагора и его школы, а также Никомаха (\cite{G}--\cite{To}).  Источником для введения средних для древнегреческих учёных стали пропорции, учение о которых активно развивалось, так как понятие пропорции использовалось тогда не только в математике, но  также и в философии, скульптуре, архитектуре и живописи (\cite{G}--\cite{To}). Например, если найти величину $x$ из нижеследующих пропорций, то как раз последовательно получим основные средние
$$
\frac{a-x}{x-b}=\frac{a}{a}=1 \Longrightarrow x=A(a,b),
$$
$$
\frac{a-x}{x-b}=\frac{a}{x} \Longrightarrow x=G(a,b),\ \
\frac{a-x}{x-b}=\frac{a}{b} \Longrightarrow x=H(a,b).
$$
Позднее были введены \textsl{среднее квадратичное}
$$
Q(a,b)=\sqrt{\frac{a^2+b^2}{2}}
$$
 и общие \textsl{степенные средние порядка $r$}, которые для набора неотрицательных (при $r \geq 0$) или  положительных (при $r < 0$) чисел определяются по формулам
\begin{equation} \label{pm}
M_r(a_1,a_2,\dots,a_n)=\left(\frac{a_1^r+a_2^r+\dots a_n^r}{n}\right)^{\frac{1}{r}},
\end{equation}
$$
M_0(a_1,a_2,\dots,a_n)=\lim_{r\to 0} M_r(a_1,a_2,\dots,a_n)=\sqrt[n]{a_1\cdot a_2\cdot\  \dots\  \cdot a_n}\ \ .
$$
При этом  справедливы формулы, выражающие основные средние через степенные. Например, для двух чисел
$$
M_{-1}(a,b)=H(a,b),\   M_0(a,b)=G(a,b),\   M_1(a,b)=A(a,b),\  M_2(a,b)=Q(a,b).
$$

Важнейшим свойством степенных средних является тот факт, что семейство этих величин монотонно возрастает по параметру $r$, или, как говорят, образует шкалу (\cite{HLP}--\cite{B}). В частности,
$$
H(a,b) \leq  G(a,b) \leq A(a,b).
$$
Эта самая известная из целого ряда теорем, в которой сравниваются между собой различные средние.

Наряду со степенными  в различных задачах потребовалось использовать и другие типы средних (\cite{B}--\cite{Sit2}). В 1938 году итальянский математик и экономист Коррадо Джини  (Corrado Gini, 1884--1965) ввёл новое семейство средних, которые впоследствии были названы его именем. Эти средние величины определяются так. Для любого набора  положительных  чисел $a=(a_1,a_2,\dots,a_n)$ введём степенные суммы
\begin{equation}\label{sum}
S_p(a)=\sum_{i=1}^n a_i^p.
\end{equation}
Тогда \textsl{средние Джини} $Gi_{p,q}(a)$, зависящие от двух вещественных параметров $p$ и $q$, определяются по формулам \begin{equation}\label{Gini}
Gi_{p,q}(a)=Gi_{q,p}(a)=
\left\{
\begin{array}{rl}
\left(\frac{S_p(a)}{S_q(a)}\right)^{\frac{1}{p-q}}, & \mbox{если } p\not =q \\ \\
\exp\left(\frac{\sum_{i=1}^n a_i^p \ln a_i}{\sum_{i=1}^n a_i^p}\right), & \mbox{если } p=q\not=0 \\ \\
M_0(a), & \mbox{если } p=q=0
\end{array} \right.
\end{equation}
При этом нетрудно показать, что в соотношениях (\ref{Gini}) второй случай получается из общего первого, а третий из второго при помощи предельных переходов. В частном случае, когда один из параметров равен нулю, среднее Джини сводится к степенному среднему (\ref{pm}), а именно справедливо соотношение
\begin{equation}\label{parcase}
Gi_{p,0}(a_1,a_2,\dots,a_n) = M_p(a_1,a_2,\dots,a_n).
\end{equation}
Мы далее будем рассматривать общий случай, положив для определённости  $p>q$.

Основным предметом изучения Коррадо Джини в экономической теории было имущественное
неравенство, для его описания он  ввёл знаменитый "\textsl{индекс Джини}"\,, для вычисления и анализа которого использовались в том числе и различные классы средних. Это одно из основных понятий современной социальной статистики \cite{Stat}, а сам К.~Джини --- один из общепризнанных её создателей. Его результаты также известны в области статистического анализа проблем демографии. Джини использовал и пропагандировал результаты своего соотечественника Вильфреда Парето.
Парето первым провёл конкретные расчёты, из которых следовало,
что в его время всего 20\% населения владело 80\% национального
богатства, и сделал вывод, что организованное таким образом государство ждёт быстрая и неизбежная гибель.
Просто отметим, что сейчас в России реальные показатели индексов Джини ещё хуже (для неимущих). Развивая  работы Парето, американский экономист О.~Лоренц разработал теорию оценки разницы доходов населения по группам (кривая Лоренца),
которая сейчас называется методологией Парето--Лоренца--Джини \cite{Stat}.
Отметим и третьего знаменитого итальянского математика и статистика 20 века --- Бруно де Финетти, разработавшего собственную оригинальную систему взглядов на основания теории вероятностей, автора известного в финансовой математике \textsl{парадокса де Финетти} (если капитал страховой компании конечен, то за бесконечное время она наверняка разоряется).

Несмотря на то, что
по запросу в интернете на фамилии Парето, Джини или де Финетти будет выдано несколько  тысяч ссылок, найти информацию о них самих практически невозможно.  Это неумное замалчивание есть проявление так называемой политкорректности и  связано с тем,
что все они были сторонниками итальянского фашизма, в котором, видимо, видели способ возрождения  Италии и её некоторую защиту от охваченного фашистской чумой могучего немецкого соседа. На идеях Парето воспитывался Муссолини, а Джини при нём занимал высокие посты и получал полную поддержку. Де Финетти прошёл путь от идейного сторонника фашизма до убеждённого коммуниста, читавшего лекции вьетнамским бойцам в джунглях под американскими бомбёжками. Разумеется, идеология фашизма во всех его проявлениях достойна однозначного осуждения, как и все его носители, кем бы сами они себя не считали. Однако это не повод замалчивать значительный вклад этих учёных в общечеловеческую культуру и науку, тем более что они никак не запятнали себя  непосредственным участием  в преступлениях. Понятие оптимума по Парето, средних Джини и индекса Джини, теория Колмогорова---Нагумо---Де Финетти квазиарифметических средних, названная   по именам её создателей, --- всё это важнейшие и необходимые составляющие современной математики и статистики.

Важный частный случай средних Джини получается из (\ref{Gini}) при
$q=p-1$. Эти средние были переоткрыты в 1971 г. Д.~Лемером и
имеют вид
$$
Le_p(a)=\frac{\sum_{i=1}^n a_i^p}{\sum_{i=1}^n a_i^{p-1}}.
$$
Их в советское время чуть позже Лемера также переоткрыл школьный учитель Ю.~М.~Фирсов \cite{Ka}.

\medskip

2. Покажем, что  средние  Джини так же как и степенные средние образуют шкалу по своим параметрам, для них справедлива  теорема сравнения. А именно, средние Джини увеличиваются, когда или параметр $p$, или параметр $q$, или оба эти параметра увеличиваются.

\textbf{Теорема.} Пусть среди положительных чисел $(a_1,a_2,\dots,a_n)$ есть хотя бы два неравных числа. Если выполнены условия
$$
p_1>q_1, \  p_2>q_2,  \ p_2\geq p_1, \  q_2 \geq q_1,
$$
и хотя бы одно из двух последних неравенств строгое, то справедливо неравенство
\begin{equation}\label{ineq1}
Gi_{p_1,q_1}(a) < Gi_{p_2,q_2}(a).
\end{equation}

\textbf{Доказательство.} Неравенство (\ref{ineq1}) в развёрнутом виде, записанное через степенные суммы (\ref{sum}),  выглядит так:
\begin{equation*}
\left(\frac{S_{p_1}(a)}{S_{q_1}(a)}\right)^{\frac{1}{p_1 - q_1}} < \left(\frac{S_{p_2}(a)}{S_{q_2}(a)}\right)^{\frac{1}{p_2 - q_2}}.
\end{equation*}
Оно упрощается, если от обеих частей взять логарифм:
\begin{equation}\label{ineq2}
\frac{1}{p_1 - q_1}\left[\ln S_{p_1}(a) - \ln S_{q_1}(a)\right] < \frac{1}{p_2 - q_2}\left[\ln S_{p_2}(a) - \ln S_{q_2}(a)\right].
\end{equation}
Введём функцию $f(p)=\ln S_p(a)$. Тогда неравенство (\ref{ineq2}) переписывается так:
\begin{equation}\label{ineq3}
\frac{f(p_1)-f(q_1)}{p_1 - q_1} < \frac{f(p_2)-f(q_2)}{p_2 - q_2}.
\end{equation}
Но (\ref{ineq3}) заведомо выполняется, если функция $f(p)$ строго выпукла. Действительно, левая часть (\ref{ineq3}) --- это тангенс наклона секущей, проходящей через точки $(p_1, f(p_1)),(q_1, f(q_1))$,
а правая часть --- это тангенс наклона второй секущей. По условиям теоремы первый тангенс меньше второго тангенса, см. рисунок 1.

\begin{center}
\includegraphics[scale=0.3,angle=-90]{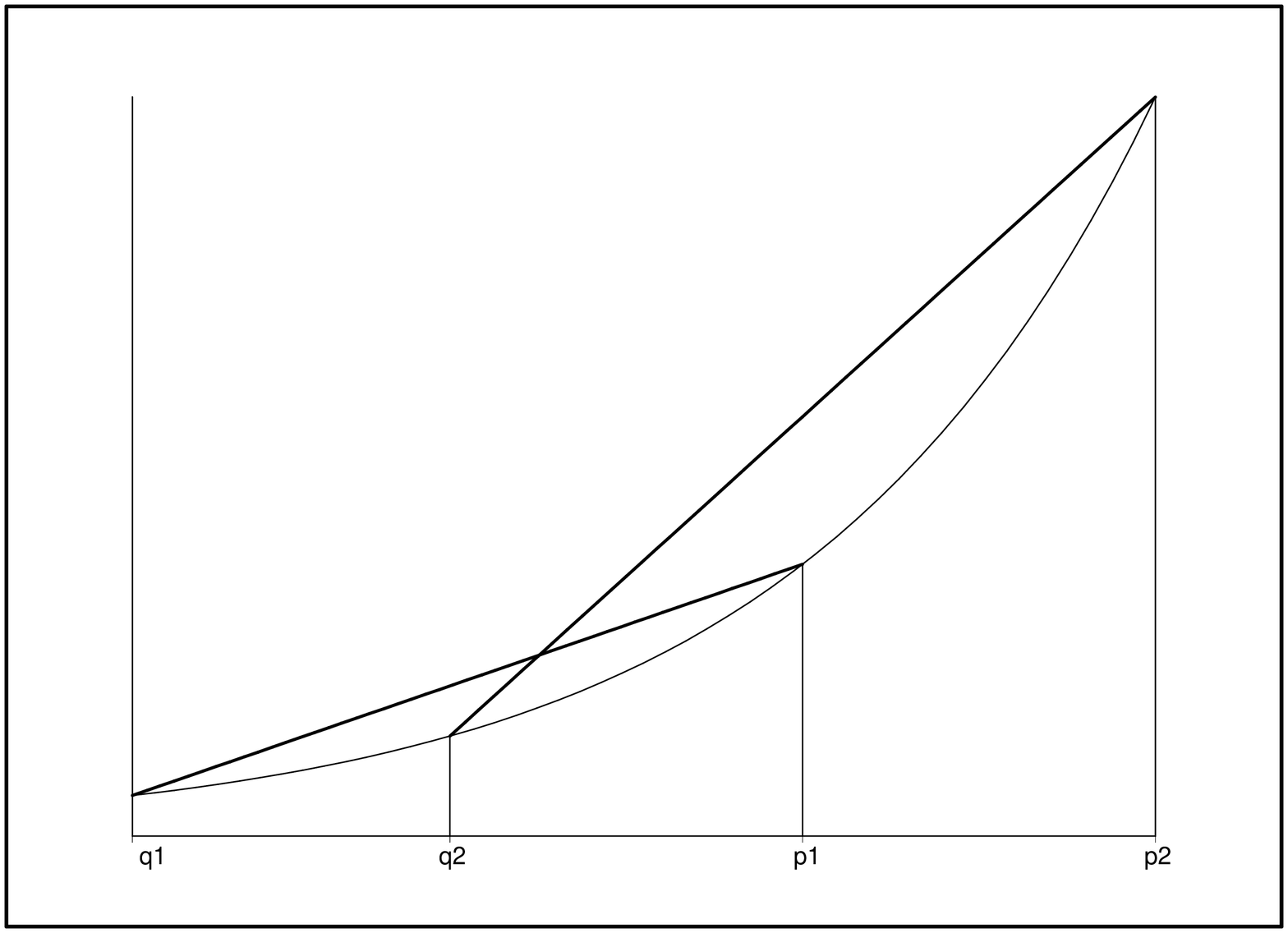}
\end{center}

\begin{center}
Рис. 1.
\end{center}

Осталось установить строгую выпуклость функции $f(p)$, для чего достаточно установить неравенство $f''(p) > 0, -\infty <p< \infty.$ Вычисляем
$$
f''(p)=\frac{1}{S^2(p)}\left[S''(p) S(p) - \left(S'(p)\right)^2\right].
$$
Это выражение будет положительным, если выполнено условие
\begin{equation}\label{log}
\left(S'(p)\right)^2 < S''(p) S(p).
\end{equation}
В развёрнутом виде  неравенство (\ref{log}) эквивалентно следующему:
$$
\left(\sum_{i=1}^n a_i^p \ln a_i\right)^2 < \sum_{i=1}^n a_i^p \sum_{i=1}^n a_i^p (\ln a_i)^2.
$$
Но последнее соотношение следует из неравенства Коши --- Буняковского для векторов с координатами $\sqrt{a_i^p}$ и  $\sqrt{a_i^p} \ln a_i$.

Покажем, что действительно неравенство (\ref{log}) всегда является строгим, хотя в использованном нами неравенстве Коши -- Буняковского возможен и случай равенства. В этом случае должно существовать такое число $\lambda$, что
$$
\sqrt{a_i^p} \ln a_i = \lambda \sqrt{a_i^p}, i=1,2,\dots, n,
$$
а это возможно только, если все числа $a_i$ равны, что противоречит условию теоремы. Следовательно, в (\ref{log}) всегда будет строгое неравенство. Поэтому при сформулированных в условиях теоремы предположениях также строгими будут неравенства (\ref{ineq3}) и (\ref{ineq1}).

Это завершает доказательство теоремы.

По мнению авторов, приведённое доказательство  проще, чем в монографии Буллена (\cite{B}, стр. 249), в том числе и потому, что не содержит внешних ссылок на  результаты других работ.

В качестве следствия рассмотрим неравенства между средними Джини и более известными степенными средними, которые получаются с учётом соотношения (\ref{parcase}).

\textbf{Следствие.} Пусть среди положительных чисел $(a_1,a_2,\dots,a_n)$ есть хотя бы два неравных числа. Тогда, если выполнены условия
$$
r \geq p > q, \  r > 0 > q,
$$
то справедливо неравенство
\begin{equation*}
Gi_{p,q}(a_1,a_2,\dots,a_n) < M_r(a_1,a_2,\dots,a_n).
\end{equation*}
Если выполнены условия
$$
p \geq r > 0, \  p > q > 0,
$$
то справедливо неравенство
\begin{equation*}
M_r(a_1,a_2,\dots,a_n) < Gi_{p,q}(a_1,a_2,\dots,a_n).
\end{equation*}

В качестве упражнения предлагаем читателю самостоятельно доказать, используя определение (\ref{Gini}), что выполнены соотношения
$$
\lim_{p\to \infty} Gi_{p,q}(a)=\max_{1\le i \le n} a_i,
\lim_{q\to -\infty} Gi_{p,q}(a)=\min_{1\le i \le n} a_i,
$$
из которых следует, что среднее Джини (с учётом доказанной теоремы о монотонности по параметрам!) находится между минимальным и максимальным из чисел $a_i$:
$$
\min_{1\le i \le n} a_i \leq Gi_{p,q}(a) \leq \max_{1\le i \le n} a_i,
$$
и таким образом является настоящим средним.

\medskip

3. Средние значения достаточно широко используются в химии, в частности, в теории полимеризации. Укажем на  использования различных средних в этом разделе химии на примере известных монографий \cite{BVE}--\cite{Fre}. Так, в монографии \cite{BVE} в главе о молекулярно--массовом распределении изучаются  различные способы усреднения по массе молекул и по числу частиц. Вводится и изучается весовое среднее арифметическое (\cite{BVE}, c.~127), для средневязкостной молекулярной массы вводится весовое степенное среднее вида (\ref{pm}), см. (\cite{BVE}, c.~130). При этом с математической точки зрения без доказательства используется монотонность весовых степенных средних по параметру (свойство шкалы средних, см.  \cite{HLP}--\cite{Sit1}), а также неравенство между весовыми средними и невесовым средним арифметическим. При изучении связей механизмов полимеризации с молекулярно - массовым распределением вводятся и изучаются свойства интегральных аналогов средних Джини  $Gi_{3,2}, Gi_{2,1}$ (\cite{BVE}, с.~184--185), с их использованием определяется так называемый коэффициент полидисперсности Шульца.

В монографии \cite{Vol}  при изучении физических методов исследования макромолекул в растворах  вводится так называемое вискозиметрическое среднее, которое выражается через молекулярные веса макромолекул c помощью среднего Джини  $Gi_{s+1,1}$ (\cite{Vol}, c.~36).

Различные средние величины также широко используются в известной монографии \cite{Fre}. Так, на с.~23--24 для средних значений молекулярных весов вводятся последовательно весовое интегральное среднее произвольного порядка, вес которого равен численной функции распределения; среднее Джини, совпадающее со средним Лемера $Gi_{q,q-1}=Le_q$, которое называется q --- средним молекулярным весом; а также степенное среднее. На с.~67 при статистическом анализе механизмов фракционирования вводится среднее Джини--Лемера $Gi_{2,1}=Le_2$ . На с.~83 при изучении гидродинамических средних весов и связанных с ними критериев полидисперсности вводятся понятия среднедиффузионного и среднеседиментационного весов, которые выражаются через весовые степенные интегральные средние некоторых отрицательных порядков, а также весовое интегральное среднее Джини $Gi_{1,1-b}$, где $b$ --- это некоторая характеристика, изменяющаяся в пределах $0<b<1$. При этом без строгого доказательства с необычной для математика формулировкой "обычно это верно"\  утверждается, что указанное среднее Джини всегда больше среднего арифметического. Отметим, что строгое доказательство этого неравенства вытекает из приведённого выше следствия при выборе параметров $p=1, q=1-b, 0<q<1, r=1$. На с.~83--84 вводится некоторая величина, связанная с молекулярными весами, которая выражается через среднее Джини--Лемера $Gi_{2-b,1-b}=Le_{2-b}$, а на  с.~229 для минимальной величины эффективного параметра получено выражение через среднее Джини $Gi_{3/2,-3/2}$. Вопросы предельного поведения средних Джини в терминах отношения степенных средних рассматриваются также на с.~252 при изучении возможностей сведения процессов полимеризации к основным статистическим классам распределений.

Таким образом,  различные классы средних значений, включая средние Джини, находят важные применения в химии при рассмотрении многочисленных вопросов теории полимеризации молекул. При этом ряд используемых в химической литературе нестрогих оценок между средними могут быть строго обоснованы с использованием доказанных в настоящей статье неравенств между средними Джини и степенными средними.

\medskip

4. В заключении отметим, что так как при нашем доказательстве используется неравенство Коши--Буняковского, то можно получить некоторые усиления доказанных неравенств с использованием  уточнений неравенства\\ Коши--Буняковского, полученных в (\cite{Sit1}--\cite{Sit2}). Авторы также использовали неравенства для средних Джини при доказательстве интересного неравенства А.С.~Гаспаряна между однородными степенными формами (\cite{SP1}--\cite{SP2}).

----------------------------------------------------\\
\vskip 20mm

АВТОРЫ:\\
\\
\\
Ситник С.М.\\
Воронежский институт МВД, Воронеж, Россия.\\
Sitnik S.M., Voronezh Institute of the Ministry of Internal Affairs of Russia.\\
mathsms@yandex.ru, pochtasms@gmail.com
\\
\\
Певный А.Б.\\
Сыктывкарский госуниверситет, Сыктывкар, Россия.\\
Pevnyi A.B., Syktyvkar State University, Syktyvkar, Russia.\\
pevnyi@syktsu.ru

\end{document}